\documentclass{article}

\usepackage{bbold}
\usepackage{amsmath}
\usepackage{amssymb}
\usepackage{amscd}
\usepackage{pstricks,pst-node}

\newtheorem{defn}{Definition}
\newtheorem{lem}[defn]{Lemma}
\newtheorem{prop}[defn]{Proposition}
\newtheorem{thm}[defn]{Theorem}
\newtheorem{pf}{Proof}
\newtheorem{exmp}[defn]{Example}

\begin{document}

\title{Splitting cycles in graphs\footnote{This work was funded by NSWCDD 
through the In-house Labratory Independent Research (ILIR) program}}

\author{David A. Johannsen and David J. Marchette\\
Naval Surface Warfare Center\\Code Q33\\17320 Dahlgren
Road\\Dahlgren, VA 22448-5100}


\maketitle

\begin{abstract}
The goal of this paper is to describe a sufficient condition on cycles
in graphs for which the edge ideal is splittable.  We give an explicit
splitting function for such ideals.
\end{abstract}



\section{Algebraic preliminaries\label{background}}
We begin by recalling some definitions and theorems from commutative
algebra.  What follows in this introductory section is contained in
the paper of Eliahou and Kervaire (\cite{EK90}).  Let $I$ be a monomial
ideal in a polynomial ring whose coefficient field has characteristic
zero.  We will denote by $\mathcal{G}(I)$ a minimal set of
generators for $I$ (in fact, \emph{the} minimal set of generators).
That is, $\mathcal{G}(I)$ is the set of all monomials in $I$ that are
not proper multiples of any monomial in $I$.  We note that if $I$ is
the edge ideal of a simple graph, $G$, then the monomials corresponding
to the edges of $G$ are such a minimal set of generators.

\begin{defn}\label{splitting}A monomial ideal, $I$, is
\emph{splittable} if $I$ is the sum of two nonzero monomial ideals, $J$
and $K$, that is, $I = J + K$, such that 
\begin{itemize}
\item[(1)]$\mathcal{G}(I)$ is the disjoint union of $\mathcal{G}(J)$ and
$\mathcal{G}(K)$
\item[(2)]there is a \emph{splitting function}
\begin{align*}
\mathcal{G}(J \cap K) &\rightarrow \mathcal{G}(J) \times \mathcal{G}(K)\\
w &\mapsto (\phi(w),\psi(w))
\end{align*}
satisfying
\begin{itemize}
\item[(a)]for all $w \in \mathcal{G}(J \cap K), w =
\mbox{lcm}(\phi(w),\psi(w))$,
\item[(b)]for every subset $S \subset \mathcal{G}(J \cap K)$, both
$\mbox{lcm}(\phi(S))$ and $\mbox{lcm}(\psi(S))$ strictly divide $\mbox{lcm}(S)$.
\end{itemize}
\end{itemize}
If $J$ and $K$ satisfy both the above properties, then we say that $I =
J+K$ is a \emph{splitting} of $I$.\end{defn}

The following theorem, which relates the minimal free resolution of
each ideal in the splitting to the minimal free resolution of their sum
(see \cite{EK90}), is of central importance in what follows.

\begin{thm}[Eliahou-Kervaire]\label{Eliahou-Kervaire}Suppose $I$ is a
splittable monomial ideal with splitting $I = J+K$. Then for all $i
\ge 0$,
\begin{displaymath}
\beta_i(I) = \beta_i(J) + \beta_i(K) + \beta_{i-1}(J \cap K).
\end{displaymath}\end{thm}

\section{Splitting cycles\label{splittingcycles}}
We now want to apply the notion of splitting monomial ideals to the
specific case of edge ideals associated to simple (and undirected)
graphs.  In particular, we want to ``split along'' a particular subgraph.
To this end, we recall some standard definitions from graph theory and we
establish the notation and conventions that we will use in what
follows.  

If $G$ is a simple graph with edge set $E$ and vertex set $V
= \{x_1, \ldots, x_n\}$, then the \emph{edge ideal} of $G$ is the
(square-free quadratic) monomial ideal $\mathcal{I}(G) = \{x_i x_j |
\{x_i, x_j\} \in E\}$ contained in the polynomial ring $R = k[x_1,
\ldots, x_n]$.  We will hereafter assume that the field $k$ has
characteristic zero (so the reader may feel free to take $k$ to be any
of $\mathbb{Q}$, $\mathbb{R}$, or $\mathbb{C}$).

We now suppose that $G$ is a graph of order $n$ that contains an
(induced) $k$-cycle with no chords, where $k \ge 4$.  We will adopt the
following notation.  Let $G = (V,E)$, where $V$ are the vertices and
$E$ are the edges in the graph $G$.  Let $U = \{u_1, \ldots, u_k\}$
denote the vertices in the $k$-cycle and $W = \{w_1, \ldots, w_{n-k}\}$
denote the vertices of $G$ not contained in the cycle.  Then, $V = U
\sqcup W$.  Now, let $\mbox{$E_U = \{u_1 u_2, u_2 u_3, \ldots, u_{k-1}
u_k, u_k u_1\}$}$; i.e., $E_U$ is the $k$-cycle contained in $G$.  Also
denote $\mbox{$E_W = \{w_p w_q ~|~ w_p, w_q \in W \mbox{ and } w_p w_q
\in E\}$}$.  Finally, let $\mbox{$E_X = \{u_i w_p ~|~ u_i \in U \mbox{
and } w_p \in W \mbox{ and } u_i w_p \in E\}$}$.  Just as we wrote for
the vertices of $G$, we can write the edges as a disjoint union, $E =
E_U \sqcup E_W \sqcup E_X$.  We will often denote the $k$-cycle by
$\mathcal{C}_k = (U,E_U)$.  Recall, also, that if $H \subset G$ is a
subgraph, then by the \emph{complement} of $H$ in $G$ we shall mean the
subgraph of $G$ consisting of all vertices of $G$ and those edges in
$G$ that are not in $H$.  Finally, to simplify the notation, we will
adopt the convention in what follows that \emph{any subscript appearing
on a ``u'' labeled vertex (i.e., a vertex in $U$) will be understood to
be $(\mbox{mod }k)$}.

We wish to parallel the work that H\`{a} and Van Tuyl did on splitting
edges in a graph (\cite{HT07}).  In what follows, we prove a sufficient
condition for a cycle to be splitting.  We begin by making the obvious
definition.

\begin{defn}We denote by $G\setminus\mathcal{C}_k$ the complement of
the $k$-cycle in $G$ (i.e., $G\setminus\mathcal{C}_k = (V, E_W \cup E_X)$).
Then, a $k$-cycle, $\mathcal{C}_k$, is a \emph{splitting cycle} if
$\mathcal{I}(G) = \mathcal{I}(\mathcal{C}_k) +
\mathcal{I}(G\setminus\mathcal{C}_k)$ is a splitting.\end{defn} 

Now, let $G$ be a graph containing an induced $k$-cycle with no chords.
We denote the edge ideal of the graph, $G$, by $I := \mathcal{I}(G)$,
and we then let $J := \mathcal{I}(\mathcal{C}_k)$ and $K :=
\mathcal{I}(G\setminus\mathcal{C}_k)$.

\begin{lem}With the notation above, $I = J+K$.\end{lem}

\begin{pf}The proof of this lemma is essentially a tautology.  Note that
$\mbox{$\mathcal{G}(J) = E_U$}$ and that $\mbox{$\mathcal{G}(K) = E_W
\sqcup E_X$}$.  Thus, $\mathcal{G}(I) = \mathcal{G}(J) \sqcup
\mathcal{G}(K)$, and the result is immediate.  $\Box$\end{pf}

We now want to characterize the minimal generating set of the
intersection of the ideals $J$ and $K$, $\mathcal{G}(J \cap K)$.  
Before stating the lemma we recall some facts about the minimal
generating sets of monomial ideals.  First, if $J$ and $K$ are
monomial ideals, then the intersection, $J \cap K$ is generated by the
set 
\begin{displaymath}
\{\mbox{lcm}(f,g) ~|~ f \in J, g \in K\}
\end{displaymath}
(see, for example, \cite{V01}).  Then, one reduces this generating
set to a minimal one by checking for pairwise divisibility.  

\begin{lem}\label{notlem}Let $I, J$, and $K$ be as defined above.
Furthermore, let
\begin{eqnarray*}
A &=& \{u_i u_{i+1} w_p ~|~ u_i u_{i+1} \in E_U \mbox{ and
} (u_i w_p \in E_X \mbox{ or } u_{i+1} w_p \in E_X)\}\\
B &=& \{u_i u_{i+1} u_j w_p ~|~ u_i u_{i+1} \in E_U \mbox{
and } u_j w_p \in E_X \mbox{ and } u_i u_{i+1} w_p, u_{i+1} u_j w_p
\notin A\}\\
C &=& \{u_i u_{i+1} w_p w_q ~|~ u_i u_{i+1} \in E_U \mbox{ and } w_p
w_q \in E_W \mbox{ and } u_i u_{i+1} w_p, u_i u_{i+1} w_q \notin A \}
\end{eqnarray*}
Then
$\mathcal{G}(J \cap K) =  A \sqcup B \sqcup C$.\end{lem}

\begin{pf} We provide a constructive proof of this lemma.  To construct
the set $\mbox{$\mathcal{G}(J \cap K)$}$ one begins by considering all possible
products of elements of $\mathcal{G}(J)$ and $\mathcal{G}(K)$.  From
this set, remove any duplicate monomials.  Then, using the
characterization of \emph{the} minimal generating set as the maximal
set such that no element divides any other, we reduce this collection.

We begin by identifying the monomials in our proposed $\mathcal{G}(J
\cap K)$ that contain a square (i.e., a $u_i^2$ for some $i$).  We
replace these degree four monomials with the degree three monomial in
which $u_i$ replaces the $u_i^2$.  Note that there can be no terms
containing a $w_p^2$ since $\mathcal{G}(J)$ contains only $u$ variables
and our graph is assumed to be simple (so there is no $p$ such that
$w_p^2 \in \mathcal{G}(K)$). Then the degree three monomials are
exactly the set $A$.

Continuing, we next eliminate the redundant degree four monomials
that are of degree three in the $u$ variables and degree one in the $w$
variables.  These are precisely the the degree four monomials that are
contained in the ideal generated by $A$ (i.e., that are multiples of
some element of $A$).  So, we want to remove redundant monomials of the
form $u_i u_{i+1} u_j w_p$; that is, we want to determine when $u_i
u_{i+1} u_j w_p$ is a multiple of some $u_k u_{k+1} w_q \in A$, some
$k$.  First, note that we must have $p = q$.  Suppose not; i.e.,
suppose that $w_p \neq w_q$.  We are assuming that $u_k u_{k+1} w_q
\mid u_i u_{i+1} u_j w_p$.  Now $U \cap W = \emptyset$ implies $u_k
\nmid w_p$ (and similarly $u_{k+1} \nmid w_p$, of course) and, by
assumption, $w_q \nmid w_p$.  Thus, $u_k u_{k+1} w_q \mid u_i u_{i+1}
u_j$.  However, both monomials are of degree three and $w_q \nmid u_m$
for any $m$.  This is a contradiction.  Thus, $p = q$.  Now, $u_k
u_{k+1} w_p \mid u_i u_{i+1} u_j w_p$ if and only if $u_k u_{k+1} \mid
u_i u_{i+1} u_j$.  Thus, either $k = i$, or $k = i+1$ and $j = k+1
(=i+2)$.  This is precisely the set $B$.

The final step in producing $\mathcal{G}(J \cap K)$ is to eliminate the
degree four monomials that are of degree two in both $u$ and $w$ and
that are contained in the ideal generated by $A$ (note, we obviously
needn't consider multiples of elements of the set $B$).  Elements of
degree two in both $u$ and $w$ must be of the form $u_i u_{i+1} w_p
w_q$, as they are products of elements of $\mathcal{G}(J)$ and
$\mathcal{G}(K)$.  Thus, an element of this form is a multiple of an
element of $A$ if $u_i u_{i+1} w_p \in A$ or $u_i u_{i+1} w_q \in A$.
This is precisely the set $C$.  

The very nature of our construction (i.e., eliminating
redundant monomials) guarantees that the decomposition is as a disjoint
union.$\Box$\end{pf}

We are now ready to state and prove the main proposition about
splitting cycles.  In the proposition immediately below, we give a
sufficient condition on a cycle for the cycle to be splitting.

\begin{prop}\label{mainprop}If there is no $i \in \{1,\ldots,k\}$ such
that $u_i w_p, u_{i+1}w_q \in \mbox{E}$ then $\mathcal{C}_k$ is a
splitting cycle.  Equivalently, the cycle, $\mathcal{C}_k$, is
splitting if it contains no adjacent vertices of degree greater than
2.\end{prop}

\begin{pf}We first comment that the equivalence of the two statements
of the proposition is immediate and obvious. So, we now proceed with
the proof of the first statement in the proposition.

Using the notation introduced in Lemma \ref{notlem} we first define a
splitting function, $\mbox{$\mathcal{G}(J \cap K) \rightarrow
\mathcal{G}(J) \times \mathcal{G}(K)$}$, by
\begin{displaymath}
w \mapsto (\phi(w), \psi(w)) = \left\{ \begin{array}{c@{\quad,\quad}l}
(u_i u_{i+1},u_i w_p) & w \in A \mbox{ and }  u_i w_p \in E_X \\
(u_i u_{i+1},u_{i+1} w_p) & w \in A \mbox{ and }  u_{i+1} w_p \in E_X \\
(u_i u_{i+1}, u_j w_p) & w \in B \\
(u_i u_{i+1}, w_p w_q) & w \in C.
\end{array} \right.
\end{displaymath}
Now, we need to verify that our function satisfies conditions (a) and
(b) given in Definition \ref{splitting}.  Note the condition (a) is
an immediate consequence of the decomposition of $\mathcal{G}(J \cap K)$
given in Lemma \ref{notlem} above.  So, let $S \subset \mathcal{G}(J
\cap K)$.  Our description of $\mathcal{G}(J \cap K)$ in Lemma
\ref{notlem} implies that each element of $S$ is divisible by some
$w_p$ while no element of $\phi(S)$ is divisible by any $w_p$.  Thus,
$\mbox{\emph{lcm}}(\phi(S))$ strictly divides $\mbox{\emph{lcm}}(S)$.

We introduce one more notational convenience. We define the (open)
neighborhood of a $k$-cycle, $\mathcal{C}_k$, to be the union of the
(open) neighborhoods of the vertices in the cycle minus the vertices in
the cycle itself.  That is
\begin{displaymath}
\mathcal{N}(\mathcal{C}_k) := \left( \cup_{u_i \in \mathcal{C}_k}
\mathcal{N}(u_i) \right) \setminus U
\end{displaymath}
We can now conveniently characterize $\psi(S)$ as follows:
\begin{displaymath}
\psi(S) = \{u_j w_p ~|~ w_p \in \mathcal{N}(\mathcal{C}_k)\} \cup \{w_p
w_q ~|~ w_p, w_q \notin \mathcal{N}(\mathcal{C}_k)\}
\end{displaymath}
We now note that, by our assumption, if for some $j \in \{1, \ldots, k
\}$ we have $u_j w_p \in \psi(S)$, then we must have $u_{j-1} w_r,
u_{j+1} w_s \notin \psi(S)$ for any $r,s \in \{1, \ldots, n-k \}$.
However, for any $S \subset \mathcal{G}(J \cap K)$ we must have some $i
\in \{1, \ldots, k\}$ for which $u_i u_{i+1}$ divides
$\mbox{\emph{lcm}}(S)$. Thus, the divisibility is strict. $\Box$
\end{pf}

We should emphasize that the splitting that we have constructed in the proof
of Proposition \ref{mainprop} is valid only when the cycle has no
chords.  That is, if $\mathcal{C}_k$ has an edge $\{u_i,u_j\}$ where
$i$ and $j$ are not consecutive integers (mod ($k$)), then the
decomposition of the minimal generating set of the intersection,
$\mathcal{G}(J \cap K)$, given in Lemma \ref{notlem} does not hold.
Consequently the splitting function given in Proposition \ref{mainprop}
is not actually a splitting.  In fact, in this case we won't have even
specified a function on $\mathcal{G}(J \cap K)$, as monomials of the
form $u_h u_i u_j$ appear in $\mathcal{G}(J \cap K)$ and our definition
of the splitting function does not indicate where such terms would be
mapped.  Furthermore, the proof of strict divisibility given in
Proposition \ref{mainprop} relies the presence of $w_p$ terms for the
strict divisibility.

In Proposition \ref{mainprop} above, we showed sufficiency of the
condition on degrees of vertices to guarantee that a cycle is
splitting.  We now provide an example to show that the condition is not
also necessary.

\begin{exmp}\label{nnec}Let $G$ be the graph given below.
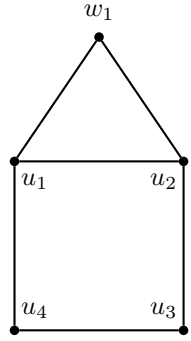
\psset{xunit=1.25in,yunit=1.25in}
\begin{center}
\begin{pspicture}(-1,-1)(1,1)
   \dotnode(.354,.354){U2}
   \dotnode(.354,-.354){U3}
   \dotnode(-.354,-.354){U4}
   \dotnode(-.354,.354){U1}
   \dotnode(0,.875){W1}
   \ncline{U1}{U2}
   \ncline{U2}{U3}
   \ncline{U3}{U4}
   \ncline{U4}{U1}
   \ncline{U1}{W1}
   \ncline{U2}{W1}
   \rput(0,.98){$w_1$}
   \rput(.27,.27){$u_2$}
   \rput(.27,-.27){$u_3$}
   \rput(-.27,-.27){$u_4$}
   \rput(-.27,.27){$u_1$}
\end{pspicture}
\end{center}
We will now show that the Eliahou-Kervaire formula holds for the Betti
numbers if we split along the four cycle.  Let $I=\;<\!u_1 u_2, u_2 u_3,
u_3 u_4, u_1 u_4, u_1 w_1, u_2 w_1\!>$.  Then, let $J = \;<\!u_1 u_2, u_2 u_3, u_3
u_4, u_1 u_4\!>$, and $K = \;<\!u_1 w_1, u_2 w_1\!>$.  It is now easy to show
that $J \cap K = \;<\!u_1 u_2 w_1, u_2 u_3 w_1, u_1 u_4 w_1\!>$.  Then we have,\\
\begin{tabular}{cccc}
$\beta(I)$ & $\beta(J)$ & $\beta(K)$ & $\beta(J \cap K)$\\
\begin{tabular}{c|ccc}
 & 1 & 2 & 3\\
\hline\\
1 & 6 & 8 & 3\\
\hline\\
total & 6 & 8 & 3
\end{tabular}
&
\begin{tabular}{c|ccc}
 & 1 & 2 & 3\\
\hline\\
1 & 4  & 4 & 1\\
\hline\\
total & 4 & 4 & 1
\end{tabular}
&
\begin{tabular}{c|cc}
 & 1 & 2\\
\hline\\
1 & 2 & 1\\
\hline\\
total & 2 & 1
\end{tabular}
&
\begin{tabular}{c|cc}
 & 1 & 2\\
\hline\\
1 & - & -\\
2 & 3 & 2\\
\hline\\
total & 3 & 2
\end{tabular}
\end{tabular}

One can easily see that the Betti numbers sum as indicated in Theorem
\ref{Eliahou-Kervaire}, while $G$ has adjacent vertices of degree
three (namely, $u_1$ and $u_2$).

We should say that it is probably not a surprise that such examples
exist.  It has long been known that there are ideal decompositions for
which one can show that no splitting function exists, but for which the
formula for the Betti numbers given in Theorem \ref{Eliahou-Kervaire} holds.
\end{exmp}


We now want to illustrate that the construction of a splitting that we
gave above fails for graphs not satisfying the hypotheses of our
proposition.  In particular, we consider a graph having adjacent
vertices of degree greater than two.  We then demonstrate that
there is no splitting function if one wishes to split along the cycle.
Furthermore, we show that the Betti number formula of Eliahou and
Kervaire does not hold for the edge ideal of this graph.

\begin{exmp}Consider the graph, $G$:
\psset{xunit=1.25in,yunit=1.25in}
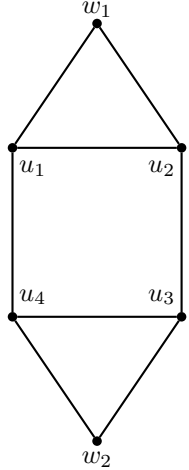
\begin{center}
\begin{pspicture}(-1,-1)(1,1)
   \dotnode(.354,.354){U1}
   \dotnode(.354,-.354){U2}
   \dotnode(-.354,-.354){U3}
   \dotnode(-.354,.354){U4}
   \dotnode(0,.875){W1}
   \dotnode(0,-.875){W2}
   \ncline{U1}{U2}
   \ncline{U2}{U3}
   \ncline{U3}{U4}
   \ncline{U4}{U1}
   \ncline{U1}{W1}
   \ncline{U4}{W1}
   \ncline{U2}{W2}
   \ncline{U3}{W2}
   \rput(.27,.27){$u_2$}
   \rput(.27,-.27){$u_3$}
   \rput(-.27,-.27){$u_4$}
   \rput(-.27,.27){$u_1$}
   \rput(0,.935){$w_1$}
   \rput(0,-.95){$w_2$}
\end{pspicture}
\end{center}
Then $\mbox{$V = \{u_1, u_2, u_3, u_4, w_1, w_2\} = \{u_1, u_2,
u_3, u_4\} \sqcup \{w_1, w_2\} = U \sqcup W$}$, and $\mbox{$E =
E_U \sqcup E_W \sqcup E_X$}$ where $\mbox{$E_U = \{u_1 u_2, u_2 u_3,
u_3 u_4, u_4 u_1\}$}$, $\mbox{$E_W = \emptyset$}$, $\mbox{$E_X =
\{u_1 w_1, u_2 w_1, u_3 w_2, u_4 w_2\}$}$.  Then the edge ideal is
given by, 
\begin{eqnarray*}
I := \mathcal{I}(G) &=& <\!u_1 u_2, u_2 u_3, u_3 u_4, u_4
u_1, u_1 w_1, u_2 w_1, u_3 w_2, u_4 w_2\!>\\
&=& \mbox{$<\!E_U \sqcup E_W \sqcup E_X\!>$}.
\end{eqnarray*}
Then, if we wish to split along the 4-cycle, we must take
\begin{align*}
J := \; <\!u_1 u_2, u_2 u_3, u_3 u_4, u_4 u_1\!> \; = \; <\!E_U\!>,\\
\intertext{and}\\
K:= \; <\! u_1 w_1, u_2 w_1, u_3 w_2, u_4 w_2\!> \; = \; <\!E_X\sqcup E_W\!>.  
\end{align*}
Then it is easy to see that $\mbox{$\mathcal{G}(J\cap K)$}= \{
u_1 u_2 w_1, u_2 u_3 w_1, u_2 u_3 w_2, u_1 u_4 w_1,  u_1 u_4 w_2\}$.

We now will show that there exists no splitting function for this
decomposition of the edge ideal (i.e., $I = J + K$ and $\mathcal{G}(I)
= \mathcal{G}(J) \sqcup \mathcal{G}(K)$).  We begin trying to construct
a function $\mathcal{G}(J\cap K) \rightarrow \mathcal{G}(J) \times
\mathcal{G}(K)$ satsifying property \emph{(a)} of Definition
\ref{splitting}.  It's quite easy to see that this requirement forces
us to define

\begin{eqnarray*}
\mathcal{G}(J\cap K) &\rightarrow& \mathcal{G}(J) \times
\mathcal{G}(K)\\
u_1 u_2 w_1 &\mapsto& \left(u_1 u_2, \begin{array}{c}u_1
w_1\\\mbox{or}\\u_2 w_1\end{array}\right)\\
u_1 u_4 w_1 &\mapsto& (u_1 u_4, u_1 w_1)\\
u_2 u_3 w_1 &\mapsto& (u_2 u_3, u_2 w_1)\\
u_1 u_4 w_2 &\mapsto& (u_1 u_4, u_4 w_2)\\
u_2 u_3 w_2 &\mapsto& (u_2 u_3, u_3 w_2)\\
\end{eqnarray*}

Now, choose $S = \mathcal{G}(J\cap K)$ and we have 
\begin{eqnarray*}
\mbox{lcm}(\psi(S)) &=& u_1 u_2 u_3 u_4 w_1 w_2\\
&=& \mbox{lcm}(S).
\end{eqnarray*}

Note that the least common multiple that we computed is the same no
matter what choice we make for the image $\psi(u_1 u_2 w_1)$.  Thus,
there is no splitting function if one decomposes the edge ideal along
the 4-cycle.

This example also illustrates why the hpotheses of Proposition
\ref{mainprop} are necessary.  That is, it is precisely the requirement
that adjacent vertices of the cycle are not both connected to the
complement of the cycle that guarantees that the property \emph{(b)} of
Definition \ref{splitting} holds.

We now also show that, unlike Example \ref{nnec}, the Betti numbers do
not sum as in the Eliahou-Kervaire formula.  Then, one computes the
Betti numbers of $I$, $J$, and $K$ to be:\\
\begin{center}
\begin{tabular}{ccc}
$\beta(I)$ & $\beta(J)$ & $\beta(K)$\\
\begin{tabular}{c|cccc}
 & 1 & 2 & 3 & 4 \\
\hline\\
1 & 8 & 12 & 5 & - \\
2 & - &  2 & 4 & 2 \\
\hline\\
total & 8 & 14 & 9 & 2 
\end{tabular}
&
\begin{tabular}{c|ccc}
 & 1 & 2 & 3\\
\hline\\
1 & 4 & 4 & 1\\
\hline\\
total & 4 & 4 & 1
\end{tabular}
&
\begin{tabular}{c|cccc}
 & 1 & 2 & 3 & 4\\
\hline\\
1 & 4 & 2 & - & -\\
2 & - & 4 & 4 & 1\\
\hline\\
total & 4 & 6 & 4 & 1
\end{tabular}
\end{tabular}
\end{center}

The Betti numbers of $(J \cap K)$ are:\\
\begin{center}
\begin{tabular}{c|ccc}
 & 1  & 2 & 3 \\
\hline\\
1 & - & - & - \\
2 & 6 & 6 & - \\
3 & - & - & 1 \\
\hline\\
total & 6 & 6 & 1
\end{tabular}
\end{center}

Now, one can easily see that, for $\beta_4$, the Eliahou-Kervaire
formula holds: 
\begin{eqnarray*}
\beta_4 (I) = \beta_4 (J) + \beta_4 (K) + \beta_3 (J \cap K),
\end{eqnarray*}
However, the formula does not hold otherwise:
\begin{eqnarray*}
\beta_3 (I) \not= \beta_3 (J) + \beta_3 (K) + \beta_2 (J \cap K)\\  
\beta_2 (I) \not= \beta_2 (J) + \beta_2 (K) + \beta_1 (J \cap K).
\end{eqnarray*}
Thus, the cycle is neither splitting nor does the sum formula for Betti
numbers hold.
\end{exmp}

Thus, we have shown a sufficient condition for a cycle to be splitting, and
provided an example showing that our condition is not also necessary.
If one hopes to use our result (together with the result of H\`{a} and
Van Tuyl, \cite{HT07}) to provide a recursive algorithm for the
computation of the Betti numbers of graphs that are not chordal, but
that contain a cycle that satisfies the hypotheses of Proposition
\ref{mainprop}, then one must be able to compute $\beta_{i-1}(J \cap
K)$ for these cycles.  This is the problem to which we are currently
turning our attention.

Though a recursive algorithm to compute the Betti numbers of graphs
that are not chordal is still out of reach, as an application of our
proposition we can now give the Betti numbers for a special class of
graphs that are not chordal.  In the following example, we derive the
formula for the Betti numbers of a wheel graph with an odd number of
vertices and every other spoke missing.

\begin{exmp}
As a non-trivial example, consider the graph
defined as follows: $G=(V,E)$ with $V=\{w,u_1,\ldots,u_{2k}\}$,
$E=\{u_1u_2,\ldots,u_{2k-1}u_{2k},u_{2k}u_1,wu_2,wu_4, \ldots,
wu_{2k}\}$.
This is a wheel graph with an odd number of vertices and every other
spoke
missing. In this case, $J$ corresponds to $C_{2k}$ and $K$ to
$S_k=\mbox{star}(k)$, the subgraph consisting of the
hub and remaining $k$ spokes of the wheel.
It is easy to see that $J\cap K = w \mathcal{I}(C_{2k})$.
Thus $\beta_{i,j}(G) = \beta_{i,j}(C_{2k})
+ \beta_{i,j}(S_k) + \beta_{i-1,j}(C_{2k})$.
Formulas for each of these is known, see \cite{J}.
Putting this together with the recursive algorithm of \cite{HT07},
we see that for $k>1$, the wheel graph $W(n)$ on $n=2k+1$ vertices
has Betti numbers given by
$$
\beta_{i,j}(W) = \left\{\begin{array}{lc}
\beta_{i,j}(C_{2k})+ \beta_{i,j}(S_k) +
\beta_{i-1,j}(C_{2k})+\sum\limits_{a=k}^{2k-1}{a\choose i},&j=i+2\\
\beta_{i,j}(C_{2k})+ \beta_{i,j}(S_k) +
\beta_{i-1,j}(C_{2k}),&\mbox{else}
\end{array}
\right.
$$

\end{exmp}


\begin{thebibliography}{}

\bibitem{EK90} S. Eliahou and M. Kervaire, \textsl{Minimal resolutions
of some monomial ideals}, J. Algebra \textbf{129} (1990), no. 1, 1-25.
MR 1037391 (91b:13019)

\bibitem{HT07} H.T. H\`{a} and A. Van Tuyl, \textsl{Splittable ideals and
the resolutions of monomial ideals}, J. Algebra \textbf{309} (2007),
no. 1, 405-425. MR 2301246 (2008a:13016)

\bibitem{V01} R. Villarreal, \textsl{Monomial Algebras}, Marcel Dekker,
New York, 2001.

\bibitem{J} S. Jacques, \textsl{Betti numbers of graph ideals}, University of
Sheffield, Ph.D. Thesis. math arXiv: 0410107v1.
\end{thebibliography}
\end{document}